\documentclass[preprint,12pt]{elsarticle}

\usepackage{amssymb, bm}
\usepackage{amsthm, amsmath}

\usepackage{graphicx,color}
\usepackage{caption}
\usepackage{subcaption}
\usepackage{multirow}
\usepackage{mathtools,xparse}
\usepackage{algorithm}
\usepackage{algorithmic}

\usepackage{listings}
\usepackage[numbers]{natbib}

\newtheorem{thm}{Theorem}

\makeatletter
\newenvironment{smallermatrix}[1][c]
{\null\,\vcenter\bgroup
  \Let@\restore@math@cr\default@tag
  \baselineskip0pt \lineskip0.4pt \lineskiplimit0pt
  \ialign\bgroup\if#1l\else\hfil\fi$\m@th\scriptstyle##$\if#1r\else\hfil\fi&&\thickspace\hfil
  $\m@th\scriptstyle##$\hfil\crcr
}{%
  \crcr\egroup\egroup\,%
}
\makeatother

\NewDocumentCommand{\ts}{O{c} e{^?_}}{
  \begin{smallermatrix}[#1]
  \mathstrut\IfValueT{#2}{#2} \\
  \mathstrut\IfValueT{#3}{#3} \\
  \mathstrut\IfValueT{#4}{#4}
  \end{smallermatrix}%
}

\journal{Arxiv}

\begin{document}

\begin{frontmatter}



\title{First-order empirical interpolation method for real-time solution of parametric time-dependent nonlinear PDEs}




\author[inst2]{Ngoc Cuong Nguyen}

\affiliation[inst2]{organization={Center for Computational Engineering, Department of Aeronautics and Astronautics, Massachusetts Institute of Technology},
            addressline={77 Massachusetts
Avenue}, 
            city={Cambridge},
            state={MA},
            postcode={02139}, 
            country={USA}}
            

\begin{abstract}
We present a model reduction approach for the real-time solution of time-dependent nonlinear partial differential equations (PDEs) with parametric dependencies. The approach integrates several ingredients to develop efficient and accurate reduced-order models. Proper orthogonal decomposition is used to construct a reduced-basis (RB) space which provides a rapidly convergent approximation of the parametric solution manifold. The Galerkin projection is employed to reduce the dimensionality of the problem by projecting the weak formulation of the governing PDEs onto the RB space. A major challenge in model reduction for nonlinear PDEs is the efficient treatment of nonlinear terms, which we address by unifying the implementation of several hyperreduction methods. We introduce a first-order empirical interpolation method to approximate the nonlinear terms and recover the computational efficiency. We demonstrate the effectiveness of our methodology through its application to the Allen–Cahn equation, which models phase separation processes, and the Buckley–Leverett equation, which describes two-phase fluid flow in porous media. Numerical results highlight the accuracy, efficiency, and stability of the proposed approach. 
\end{abstract}


\end{frontmatter}


\section{Introduction}
\label{sec:intro}


One of the primary challenges in solving time-dependent nonlinear partial differential equations (PDEs) is the high computational cost associated with resolving both spatial and temporal scales while accurately capturing the nonlinear interactions. For instance, the Allen-Cahn equation involves tracking evolving phase boundaries, which requires fine spatial resolution near the interfaces, while the Buckley-Leverett equation demands high accuracy in capturing sharp shocks and fronts associated with fluid displacement. Moreover, in parametric settings, such as optimizing oil recovery strategies or predicting material behavior under different environmental conditions, a large number of simulations must be run for varying parameters, that can drastically  increase the computational burden. This can make real-time simulation and decision-making impractical without model reduction techniques that retain accuracy while drastically reducing the computational cost.

The Allen-Cahn equation is a second-order nonlinear parabolic partial differential equation (PDE) that describes the evolution of phase boundaries in multi-component systems, particularly in the context of phase separation and interfacial dynamics \cite{Allen1979}. This equation has been applied to various physical problems, including crystal growth \cite{Shah2016}, image segmentation \cite{Benes2004}, and motion by mean curvature flows \cite{Feng2003}. It has become a foundational model in the diffuse interface approach, which is used to study phase transitions and interfacial dynamics in material science.  Originally introduced by Allen and Cahn \cite{Allen1979} to model anti-phase domain coarsening in binary alloys, the equation describes the motion of interfaces between phases, which evolve to minimize interfacial energy over time. The interaction between diffusion, which smooths spatial variations, and reaction terms, which drive the system towards stable phases, leads to the formation of characteristic phase separation patterns. The interface movement is influenced by factors such as curvature and interfacial tension, making the Allen-Cahn equation  relevant for studying phenomena like grain boundary motion \cite{Cahn1997,Deckelnick2001}, phase transitions \cite{Katsoulakis1995}, and other processes where interfacial dynamics play a key role.

The Buckley-Leverett equation is a second-order nonlinear convection-diffusion equation that models two-phase flow in porous media, particularly in the context of oil recovery and hydrology \cite{Buckley1942}. It describes the displacement of one immiscible fluid (typically water) by another (typically oil or gas) in a porous medium, accounting for saturation levels of the fluids over time and space. The equation arises from the conservation of mass for immiscible fluids and incorporates the effects of relative permeability and capillary pressure to capture the dynamics of the flow. A key feature of the Buckley-Leverett equation is the formation of sharp fronts or shocks, representing the interface between the displacing and displaced fluids, with the velocity of the shock dependent on the saturation levels and flow properties. Despite its simplicity, the Buckley-Leverett equation captures the essential physics of two-phase flow in porous media, making it a useful tool in studies of subsurface fluid dynamics, including two-phase flow in porous media \cite{Spayd2011,Cogswell2017}, oil recovery \cite{Welge1952}, groundwater contamination \cite{Corapcioglu1990} and CO$_2$ sequestration \cite{Hayek2009}.

Efficient and accurate numerical methods for solving the Allen-Cahn and Buckley-Leverett equations are of significant interest because of their broad applications.  Both equations, though different in their mathematical structure and physical applications, share common features that make their computational treatment challenging. These include the presence of sharp interfaces or fronts, nonlinearity, and the need to capture complex dynamics over a range of spatial and temporal scales. Traditional numerical methods, while capable of providing accurate solutions, often become computationally expensive when applied to large-scale or real-time problems with parametric variations. To make real-time solutions of the Allen-Cahn and Buckley-Leverett equations feasible in the context of many-query studies, reduced-order modeling is crucial.


Reduced-basis (RB) methods have been widely used to achieve rapid and accurate solutions of parametrized PDEs \cite{Chen2009, Chen2010a, deparis07, EftangPateraRonquist10,Eftang2012a, Huynh2007a, Huynh2010, PhuongHuynh2013, Karcher2018, Knezevic2011, Knezevic2010,  nguyen04:_handb_mater_model, Nguyen2007, Nguyen_SantaFE08, Calcolo, ARCME, Rozza05_apnum, Sen2006b, veroy04:_certif_navier_stokes, Veroy2002, Vidal-Codina2014, Vidal-Codina2018a}. However, nonlinear PDEs like the Allen-Cahn and Buckley-Leverett equations introduce a challenge for reduced-order modeling. Nonlinear terms pose difficulties because  the high computational cost associated with evaluating the nonlinear terms often undermines the efficiency gain of model reduction. To address this, hyperreduction techniques have been developed to provide efficient treatment of the nonlinear terms. Hyperreduction methods, such as empirical interpolation methods (EIM) \cite{Barrault2004,Grepl2007,Nguyen2007,Drohmann2012,Manzoni2012,Hesthaven2014,Hesthaven2022,Chen2021,Eftang2010b}, best-points interpolation method \cite{Nguyen2008a,Nguyen2008d}, generalized empirical interpolation method \cite{Maday2015a,Maday2013},  empirical quadrature methods (EQM) \cite{Patera2017,Hernandez2017,Yano2019a,Taddei2021,Sleeman2022,Mirhoseini2024}, energy-conserving sampling and weighting \cite{Farhat2014, Farhat2015}, gappy-POD \cite{everson95karhunenloeve,Galbally2010}, and integral interpolation methods \cite{Carlberg2011,Chaturantabut2011,Drohmann2012, Kerfriden2013,Radermacher2016}, can be used to approximate the nonlinear terms with low computational cost, making real-time solution of parametrized nonlinear PDEs  feasible.

 
In this work, we propose a model reduction approach that integrates POD, Galerkin projection, and hyperreduction techniques for the efficient and accurate solution of time-dependent nonlinear PDEs with parametric dependencies. We introduce a first-order empirical interpolation method (FOEIM) that approximates the nonlinear terms while retaining the computational efficiency of the reduced-order model. This approach can handle the complexities of nonlinear PDEs while achieving substantial computational speedups.  We consider the Allen-Cahn and Buckley-Leverett equations to demonstrate the effectiveness of the proposed approach across different scenarios. Through numerical experiments, we show that it delivers accurate solutions at speeds several orders of magnitude faster than traditional numerical methods. 


Recent advancements in model reduction for time-dependent nonlinear PDEs have seen significant progress in several areas. One prominent direction is model reduction on nonlinear manifolds, which originated with the work \cite{Rutzmoser2017} on reduced-order modeling of nonlinear structural dynamics using quadratic manifolds. This concept was further extended through the use of deep convolutional autoencoders by Lee and Carlberg \cite{Lee2020}. Recently, quadratic manifolds are further developed to address the challenges posed by the Kolmogorov barrier in  model order reduction \cite{Geelen2023,Barnett2022}. Another approach is model reduction through lifting or variable transformations, as demonstrated by Kramer and Willcox \cite{Kramer2019}, where nonlinear systems are reformulated in a quadratic framework, making the reduced model more tractable. This method was expanded in the "Lift \& Learn" framework \cite{Qian2020} and operator inference techniques \cite{McQuarrie2023,Kramer2024}, leveraging physics-informed machine learning for large-scale nonlinear dynamical systems.

The paper is organized as follows. In Section 2, we present a unified framework to implement various hyper-reduction  methods for time-dependent nonlinear PDEs with parameter dependencies. In Section 3, we introduce the first-order empirical interpolation method to provide an efficient treatment of nonlinear terms.  Numerical results are presented in Section 4 to demonstrate the accuracy, efficiency, and stability of our approach. Finally, in Section 5, we provide a number of concluding remarks on future work.

\section{Model reduction methods}

In this section, we focus on model reduction techniques for parametrized time-dependent nonlinear PDEs. The full-order model (FOM) is first constructed by using the finite element method (FEM) to solve the PDE at selected parameter values. Proper orthogonal decomposition (POD) is applied to obtain a reduced basis (RB) space from the FOM solutions. A Galerkin projection is then performed to project the FOM onto the RB space to generate a reduced-order model (ROM). Although the POD-Galerkin method reduces dimensionality significantly, it remains computationally expensive to solve the ROM system due to the non-linear terms. Hyperreduction is used to reduce the computational cost of evaluating the nonlinear terms. We present a unified implementation of several hyperreduction techniques for parametric time-dependent nonlinear PDEs.

\subsection{Finite element approximation}
\label{section3.1}

The weak formulation of a time-dependent nonlinear PDE can be stated
as follows: given any $\bm \mu \in {\cal D} \subset \mathbb{R}^P$, we
evaluate $s^{\rm e}(t,\bm \mu) = \ell^O(u^{\rm e}(t,\bm \mu);\bm \mu)$, where $u^{\rm e}(t,\bm \mu) \in
X^{\rm e}$ is the solution of
\begin{equation}
m(\dot{u}^{\rm e}(\bm \mu),v) + a(u^{\rm e}(\bm \mu),v;\bm \mu) = \ell(v; \bm \mu), \quad \forall v \in X^{\rm e},  
\label{eq:1}
\end{equation}
with initial condition $u(\bm x, 0, \bm \mu) = u_0(\bm x, \bm \mu)$. Here $\mathcal{D}$ is the parameter domain in which our $P$-tuple parameter $\bm \mu$ resides; $X^{\rm e}(\Omega)$ is an appropriate Hilbert space;
$\Omega$ is a bounded domain in $\mathbb{R}^D$ with Lipschitz continuous boundary $\partial \Omega$; $\ell(\cdot;\bm \mu), \ \ell^O(\cdot;\bm \mu)$ are $X^{\rm
  e}$-continuous linear functionals; $m(\cdot, \cdot)$ is a  symmetric positive-definite bilinear form;  and $a(\cdot,\cdot;\bm \mu)$ is a variational form of the parameterized PDE operator. We assume that $m$, $\ell$ and $\ell^O$ are independent of $\bm \mu$, and that the variational form $a$ may be expressed as
\begin{equation}
a(w,v; \bm \mu) = \sum_{q=1}^{Q} \Theta^q(\bm \mu) a^q(w,v) + b(w,v; \bm \mu)  
\end{equation} 
where $a^q(\cdot, \cdot)$ are $\bm \mu$-independent bilinear forms, and $\Theta^q(\bm \mu)$ are $\bm \mu$-dependent functions. The variational form $b$ is assumed to be 
\begin{equation}
b(w,v; \bm \mu) = \int_{\Omega} g(w,  \bm \mu) v \, d  \bm x + \int_{\Omega} \bm f(w, \bm \mu) \cdot \nabla v \, d  \bm x
\end{equation} 
where $g$ is a scalar nonlinear function of $w$, and $\bm f$ is a vector-valued nonlinear function of $w$.

Let $X \in X^{\rm e}$ be a finite element (FE) approximation space of dimension $\mathcal{N}$. We divide the time interval $(0, T]$  into $I$ subintervals of equal length $\Delta t = T/I$ and define $t_i = i \Delta t$, $0 \le i \le I$, and $\mathbb{I} = \{1, 2, \ldots, I\}$. We consider the Backward-Euler scheme for the time integration. The FE approximation $u(t_i,\bm \mu) \in X$ of the exact solution $u^{\rm e}(t_i, \bm \mu)$ is the solution of
\begin{equation}
m(u(t_i,\bm \mu),v) + \Delta t \, a(u(t_i,\bm \mu),v;\bm \mu) = \Delta t \,  \ell(v) - m(u(t_{i-1},\bm \mu),v) , 
\label{eq:7-6c}
\end{equation}
for all $v \in X$, $i \in \mathbb{I}$, and
\begin{equation}
\int_\Omega u(\bm x, 0,\bm \mu) v = \int_\Omega u_0(\bm x, \bm \mu) v,  \quad \forall  v
\in X . 
\label{eq:7-6d}
\end{equation}
The output of interest is then evaluated as $s(t_i,\bm \mu) = \ell^{O} (u(t_i,\bm \mu))$. 
We  assume that the FE discretization is sufficiently rich such that $u(t_i,\bm \mu)$ and $u^{\rm e}(t_i,\bm \mu)$
and hence $s(t_i,\bm \mu)$ and $s^{\rm e}(t_i,\bm \mu)$ are indistinguishable at the accuracy level of interest.

To construct a ROM, we compute a set of solutions $u(\bm x, t_i, \bm \mu_j), 0 \le i \le I$, at multiple parameter points $\bm \mu_j$ in a parameter sample $S_J = \{\bm \mu_j \in \mathcal{D}, 1 \le j \le J\}$. Typically, these points are selected using a greedy sampling method, which strategically explores the parameter domain to minimize the error in the ROM approximation. This process involves  error estimators that guide the selection of the most informative parameter points to improve the accuracy of the ROM. Once this set of solutions (or snapshots) is computed, POD is employed to extract the most important modes to generate a RB space that captures the dominant dynamics of the system. This RB space, in turn, forms the foundation of the ROM.

\subsection{Proper orthogonal decomposition}
\label{section2.1b}

Let $\{\zeta_k = u(t_i, \bm \mu_j)\}_{k=1}^K$ be the set of solution snapshots, where $K = IJ$ represents the total number of snapshots. The goal of POD is to generate an optimal basis set $\Phi_N = \mbox{span} \{\phi_n, 1 \le n \le N\}$ from this snapshot set. The basis set is determined by minimizing the average squared error between the snapshot set and its projection as:
\begin{equation}
\min \sum_{k=1}^{K}  \Big\|\zeta_{k} -  \frac{(\zeta_{k}, \phi)}{\|\phi\|^2} \phi \Big\|^2 \ .
\label{eq3a:8}
\end{equation} 
This minimization problem is equivalent to the maximization problem~
\begin{equation}
\max   \sum_{k=1}^{K}   |\left(\phi, \zeta_{k}\right)|^2
\label{eq3a:4}
\end{equation}
subject to the constraints $\|\phi\|^2 = 1$. It is shown in \cite{Holmes2012} (see Chapter 3) that that the problem~(\ref{eq3a:4}) is equivalent to solving the eigenfunction equation 
\begin{equation}
\frac{1}{K} \sum_{k=1}^{K}   (\zeta_{k}, \phi) \zeta_{k} = \lambda \, \phi .
\label{eq3a:5}
\end{equation}
As the optimal basis functions are given by the eigenfunctions of the eigenfunction equation, they are  called empirical eigenfunctions or POD modes. The method of snapshots~\cite{sirovich87:_turbul_dynam_coher_struc_part} expresses a typical empirical eigenfunction $\phi$ as a linear combination of the snapshots 
\begin{equation}
\phi =  \sum_{k=1}^{K}   a_{k} \zeta_{k} \ .
\label{eq3a:6}
\end{equation}
Inserting (\ref{eq3a:6}) into~(\ref{eq3a:5}), we  obtain the following eigenvalue problem
\begin{equation}
\bm C \bm a = \lambda \bm a \ ,
\label{eq3a:7}
\end{equation}
where $\bm C \in \mathbb{R}^{K \times k}$ is known as the correlation matrix with entries $C_{kk'}  = \frac{1}{K} \left(\zeta_{k},\zeta_{k'} \right), 1 \le k, k' \le K$. The eigenproblem~(\ref{eq3a:7}) can then be solved for the first $N$ eigenvalues and eigenvectors from which the POD basis functions $\phi_n, 1 \le n \le N,$ are constructed by~(\ref{eq3a:6}). These basis functions are then ordered by their associated eigenvalues, where the largest eigenvalues correspond to the most significant modes. By truncating the set to the first $N$ modes such that $\sum_{k'=N+1}^{K} \lambda_{k'}/\sum_{k=1}^{K} \lambda_{k}$ is very small, we form the RB space $\Phi_N$ to capture the dominant dynamics of the system with high accuracy.


\subsection{Galerkin-Newton method}



The POD-Galerkin formulation is obtained by a standard Galerkin projection: for any given $\bm \mu \in {\cal D}$, we evaluate $s_{N}(t_i,\bm \mu) = \ell^O (u_{N}(t_i,\bm\mu))$, where $u_{N}(t_i,\bm\mu) \in \Phi_N$ is the solution of
\begin{multline}
m(u_N(t_i,\bm \mu),v) + \Delta t \, a(u_N(t_i,\bm \mu),v;\bm \mu) = \Delta t \,  \ell(v) - m(u_N(t_{i-1},\bm \mu),v) , 
\label{eq:7-7}
\end{multline} 
for all $v \in \Phi_N$, $i \in \mathbb{I}$, and
\begin{equation}
\int_\Omega u_N(\bm x, 0,\bm \mu) v = \int_\Omega u_0(\bm x, \bm \mu) v,  \quad \forall  v
\in \Phi_N . 
\label{eq:7-7d}
\end{equation}
We now express $u_N(t_i,\bm \mu) = \sum_{n=1}^N \alpha_{N,n} (t_i,\bm \mu) \phi_n$ and choose
test functions $v = \phi_j, \ 1 \leq j \leq N$, in~(\refeq{eq:7-7}), we obtain the nonlinear algebraic system
\begin{multline}
\left(\frac{\bm M_N}{\Delta t} +  \sum_{q=1}^Q \Theta^q(\bm \mu) \bm A^q_{N} \right) \bm \alpha_N(t_i,\bm \mu) + \bm f_N(\bm \alpha_N(t_i,\bm \mu)) \\ + \bm g_N(\bm \alpha_N(t_i,\bm \mu)) =  \bm l_N + \frac{\bm M_N}{\Delta t} \bm \alpha_N(t_{i-1},\bm \mu)
\label{eq:7-8}
\end{multline} 
where, for $1 \le n, j \le N$, we have
\begin{equation}
M_{N, jn} = m(\phi_j, \phi_n), \quad A^q_{N, jn} = a^q(\phi_j, \phi_n), \quad l_{N,j} = \ell(\phi_j),
\label{eq:7-9}
\end{equation} 
and 
\begin{equation}
\label{eq:7-10}
\begin{split}
f_{N, j}(\bm \alpha_N(t_i,\bm \mu)) & = \int_{\Omega} \bm f (u_{N}(t_i,\bm \mu), \bm \mu) \cdot \nabla \phi_j , \\
g_{N, j}(\bm \alpha_N(t_i,\bm \mu)) & = \int_{\Omega} g(u_{N}(t_i,\bm \mu), \bm \mu) \phi_j  .
\end{split}
\end{equation} 
Note that $\bm M_N$, $\bm A_N^q, 1 \le q \le Q,$ and $\bm l_N$ can be pre-computed and stored in the offline stage due to their independence of $\bm \mu$, whereas $\bm f_N$ and $\bm g_N$ can not be pre-computed due to the presence of the nonlinear terms.

We use Newton's method to linearize (\ref{eq:7-8}) at a current iterate $\bar{\bm \alpha}_N(t_i,\bm \mu)$ to arrive at the following linear system 
\begin{multline}
\left( \frac{\bm M_N}{\Delta t} \ +  \sum_{q=1}^Q \Theta^q(\bm \mu) \bm A^q_{N} + \bm D_N(\bar{\bm \alpha}_N(t_i,\bm \mu)) + \bm E_N(\bar{\bm \alpha}_N(t_i,\bm \mu))  \right) \delta \bm \alpha_N(t_i,\bm \mu) \\
= \bm l_N + \frac{\bm M_N}{\Delta t} \bm \alpha_N(t_{i-1},\bm \mu) - \left(\frac{\bm M_N}{\Delta t} +  \sum_{q=1}^Q \Theta^q(\bm \mu) \bm A^q_{N} \right) \bar{\bm \alpha}_N(t_i,\bm \mu) \\ - \bm f_N(\bar{\bm \alpha}_N(t_i,\bm \mu)) - \bm g_N(\bar{\bm \alpha}_N(t_i,\bm \mu))
\label{eq:7-11}
\end{multline} 
where, for $1 \le n, j \le N$, we have
\begin{equation}
\begin{split}    
D_{N, j n}(\bar{\bm \alpha}_N(t_i,\bm \mu)) = \int_{\Omega} \bm f'_u(\bar{u}_{N}(t_i,\bm \mu), \bm \mu) \phi_n \cdot \nabla \phi_j \\
E_{N, j n}(\bar{\bm \alpha}_N(t_i,\bm \mu)) = \int_{\Omega} g'_u(\bar{u}_{N}(t_i,\bm \mu), \bm \mu) \phi_n \phi_j .
\end{split}
\label{eq:7-12}
\end{equation} 
Here $\bm f'_u$ and $g'_u$ are the partial derivative of $\bm f$ and $g$ with respect to $u$, respectively. Both $\bm D_N$ and $\bm E_N$ can not be pre-computed in the offline stage because they depend on the partial derivatives of the nonlinear terms. Although the linear system (\ref{eq:7-11}) is small, solving it is computationally expensive due to the high computational cost of forming $\bm f_N, \bm g_N, \bm D_N$ and $\bm E_N$. As a result, the Galerkin-Newton method does not offer a significant speedup over the FE method.

\subsection{Hyperreduced Galerkin-Newton method}

Hyperreduction techniques aim to significantly reduce the computational cost associated with evaluating nonlinear terms in ROMs. Existing hyperreduction techniques can be broadly classified as empirical quadrature methods \cite{Farhat2014, Farhat2015,Patera2017,Hernandez2017,Yano2019a}, empirical interpolation methods \cite{Barrault2004,Grepl2007,Nguyen2007,Drohmann2012,Manzoni2012,Hesthaven2014,Hesthaven2022,Chen2021,Eftang2010b}, and integral interpolation methods \cite{everson95karhunenloeve,Galbally2010,Carlberg2011,Chaturantabut2011,Drohmann2012, Kerfriden2013,Radermacher2016}. Empirical quadrature methods construct  a set of tailored quadrature points and weights that reduce the cost of evaluating the nonlinear integrals directly. These points and weights are optimized for the specific parameterized nonlinear integrals to reduce the computational cost while maintaining accuracy of the Gauss quadrature. Instead of approximating the integrals directly, empirical interpolation methods focus on approximating the nonlinear integrands. This technique uses a set of basis functions to represent the nonlinear integrands and interpolates them at selected points. Integral interpolation methods shifts the focus from interpolating the integrands to interpolating the integrals themselves. By directly approximating the integrals with a set of interpolation  indices and basis functions, integral interpolation methods can efficiently reduce the computational cost associated with evaluating the nonlinear integrals.

Despite their differences, hyperreduction techniques share a common approach in approximating nonlinear integrals. Let us consider the approximation of the nonlinear integrals in (\ref{eq:7-10}) as an example. All of these techniques approximate and replace $\bm g_N(\bm \alpha_N(t_i,\bm \mu))$ with the following vector
\begin{equation}
\widehat{\bm g}_N(\bm \alpha_N(t_i,\bm \mu)) = \bm G_{N M} \,  \bm b_M(\bm \alpha_N(t_i,\bm \mu))   
\end{equation}
where $\bm G_{N,M} \in \mathbb{R}^{N \times M}$ is  parameter-independent and thus computed in the offline stage, while $\bm b_{M} \in \mathbb{R}^{M}$ is a parameter-dependent vector with entries
\begin{equation}
 b_{M, m}(\bm \alpha_N(t_i,\bm \mu))  = g\left( \sum_{n=1}^N \alpha_n(t_i, \bm \mu) \phi_n(\bm x^g_m), \bm \mu \right), \quad 1 \le m \le M .
\end{equation}
Here $\{\bm x^g_m\}_{m=1}^M$ represents a set of $M$ quadrature points  or interpolation points for the nonlinear function $g(w, \bm \mu)$. These techniques differ in how $\bm G_{N,M}$ and  $\{\bm x^g_m\}_{m=1}^M$ are determined during the offline stage. In the online stage, the computational cost of evaluating $\widehat{\bm g}_N(\bm \alpha_N(t_i,\bm \mu))$ scales as $O(NM)$,  where $M$  typically scales linearly with $N$. This ensures a significant reduction in computational time. 

Similarly, these techniques approximate and replace $\bm f_N(\bm \alpha_N(t_i,\bm \mu))$ with the following vector
\begin{equation}
\widehat{\bm f}_N(\bm \alpha_N(t_i,\bm \mu)) = \sum_{d=1}^D \bm F^d_{N M} \,  \bm c_M^d(\bm \alpha_N(t_i,\bm \mu))   
\end{equation}
where $\bm F^d_{N,M} \in \mathbb{R}^{N \times M}$ are  parameter-independent matrices computed in the offline stage for each dimension $d$. The parameter-dependent vectors $\bm c_{M}^d \in \mathbb{R}^{M}$ are determined by
\begin{equation}
 c^d_{M, m}(\bm \alpha_N(t_i,\bm \mu))  = f^d\left( \sum_{n=1}^N \alpha_n(t_i, \bm \mu) \phi_n(\bm x^{f^d}_m),  \bm \mu \right), \quad 1 \le m \le M .
\end{equation}
Here $\{\bm x^{f^d}_m\}_{m=1}^M$ represent the quadrature points  or interpolation points for each nonlinear term $f^d(w, \bm \mu)$. For notation simplification, the same number of quadrature/interpolation points is used to approximate the nonlinear integrals. In practice, the number of points may vary depending on the complexity and degree of nonlinearity in the integrals.

By applying hyperreduction techniques to the nonlinear integrals as described above, we obtain the following hyper-reduced ROM: for any given $\bm \mu \in {\cal D}$, we evaluate 
\begin{equation}
s_{N}(t_i,\bm \mu) = ( \bm l_N^O )^T \bm \alpha_N(t_i,\bm \mu), 
\end{equation} 
where $\bm l_N^O$ has entries $l_{N,j}^O = \ell^O(\phi_j), 1 \le j \le N,$ and $\bm \alpha_N(t_i,\bm \mu) \in \mathbb{R}^N$ is the solution of 
\begin{multline}
\left(\frac{\bm M_N}{\Delta t} +  \sum_{q=1}^Q \Theta^q(\bm \mu) \bm A^q_{N} \right) \bm \alpha_N(t_i,\bm \mu) + \sum_{d=1}^D \bm F^d_{N M} \,  \bm c_M^d(\bm \alpha_N(t_i,\bm \mu))  \\ + \bm G_{N M} \,  \bm b_M(\bm \alpha_N(t_i,\bm \mu))  =  \bm l_N + \frac{\bm M_N}{\Delta t} \bm \alpha_N(t_{i-1},\bm \mu), \quad i \in \mathbb{I} .
\label{eq:78w}
\end{multline} 
Since this nonlinear system is purely algebraic and small, it can be solved efficiently by using Newton's method. 

The offline and online stages of the hyperreduced Galerkin-Newton method are summarized in Algorithm 1 and Algorithm 2, respectively. The offline stage is expensive and performed once. All the quantities computed in the offline stage are independent of  $\bm \mu$.  In the online stage, the RB output $s_N(\bm \mu)$ is calculated for any $\bm \mu \in \mathcal{D}$. The computational cost of the online stage is $O(N^3 + (Q + M + DM ) N^2)$ for each Newton iteration. 
Hence, as required in the many-query or real-time contexts, the online complexity is independent of $\mathcal{N}$, which is the dimension of the FOM. Thus, we expect computational savings of several orders of magnitude relative to both the FOM and the Galerkin-Newton method described earlier.

\begin{algorithm}
\begin{algorithmic}[1]
\REQUIRE{The parameter sample set $S_J = \{\bm \mu_j, 1 \le j \le J\}$.}
\ENSURE{$\bm l_N, \bm l_N^O, \bm M, \bm A_N^q, \bm G_{NM}$, $\bm F^d_{NM}$.}
\STATE{Solve the FOM (\ref{eq:7-6c}) for each $\bm \mu_j \in S_J$ to obtain $u(t_i,\bm \mu_j), 1 \le I \le I$.}
\STATE{Construct a RB space $\Phi_N = \mbox{span} \{\phi_n, \ 1 \leq n \leq N\}$ by using POD.}
\STATE{Compute and store both $\bm G_{NM}$ and $\{\bm x_m^g \in \Omega\}_{m=1}^M$ for $g(w,  \bm \mu)$.}
\STATE{Cpmpute and store $\bm F_{NM}^d$ and $\{\bm x_m^{f^d} \in \Omega\}_{m=1}^M$ for $f^d(w,  \bm \mu)$.}
\STATE{Form and store $\bm l_N, \bm l_N^O, \bm M, \bm A_N^q$.}
\end{algorithmic}
\caption{Offline stage of the Hyperreduced Galerkin-Newton method.}
\end{algorithm}

\begin{algorithm}
\begin{algorithmic}[1]
\REQUIRE{Parameter point $\bm \mu \in \mathcal{D}$ and initial guess $\bar{\bm \alpha}_N(\bm \mu)$.}
\ENSURE{RB output $s_{N}(\bm \mu)$ and updated coefficients $\bar{\bm{\alpha}}_N(\bm \mu)$.}
\STATE{Linearize (\ref{eq:78w}) around $\bar{\bm{\alpha}}_N(\bm \mu)$.}
\STATE{Solve the resulting linear system to obtain $\delta {\bm \alpha}_N(\bm \mu)$.}
\STATE{Update $\bar{\bm \alpha}_N(\bm \mu) = \bar{\bm \alpha}_N(\bm \mu) + \delta {\bm \alpha}_N(\bm \mu)$.}
\STATE{If $\|\delta {\bm \alpha}_N(\bm \mu)\| \le \epsilon$, then calculate $s_{N}(\bm \mu) = (\bm l_N^O)^T \bar{\bm \alpha}_N(\bm \mu)$ and stop.
}
\STATE{Otherwise, go back to Step 1.}
\end{algorithmic}
\caption{Online stage of the Hyperreduced Galerkin-Newton method.}
\end{algorithm}

The unified implementation of hyperreduction methods offers a versatile framework for constructing efficient ROMs for parametrized nonlinear PDEs. By generalizing the approach to nonlinear integral approximations, this methodology accommodates a wide range of hyperreduction techniques, such as empirical quadrature, empirical interpolation, and integral interpolation. Each method reduces the computational burden of evaluating nonlinear terms, which often dominate the computational cost in ROMs. Furthermore, the implementation maintains flexibility, allowing the number of quadrature or interpolation points to be adapted based on the complexity of each nonlinear integral. This adaptive feature is critical for balancing between accuracy and computational efficiency. 

 In the following section, we discuss the empirical interpolation procedure for computing the parameter-independent matrix $\bm G_{NM}$ and the interpolation points $\{\bm x^g_m\}_{m=1}^M$. A similar approach is used to compute the parameter-independent matrices $\bm F_{NM}^d$ and the interpolation points $\{\bm x^{f^d}_m\}_{m=1}^M$. This paper focuses particularly on empirical interpolation methods, which offer a powerful way to approximate nonlinear terms without significantly increasing the dimensionality of the ROM. By constructing separate RB spaces for the nonlinear integrands, empirical interpolation reduces both the computational effort and the storage requirements while retaining the accuracy of the FOM.



\section{First-Order Empirical Interpolation Method}

The  empirical interpolation method (EIM) was first introduced in \cite{Barrault2004} for constructing basis functions and interpolation points to approximate parameter-dependent functions, and developing efficient RB approximation of non-affine PDEs. Shortly later, the empirical interpolation method was extended to develop efficient ROMs for nonlinear PDEs \cite{Grepl2007}. Since the pioneer work \cite{Barrault2004, Grepl2007}, the EIM has been widely used to construct efficient ROMs of nonaffine and nonlinear PDEs for  different applications \cite{Grepl2007,Nguyen2007,Galbally2010,Drohmann2012,Manzoni2012,Hesthaven2014,Kramer2019,Hesthaven2022,Chen2021}. Rigorous a posteriori error bounds for the empirical interpolation method is developed by Eftang et al. \cite{Eftang2010b}. Several attempts have been made to extend the EIM in diverse ways. The best-points interpolation method (BPIM) \cite{Nguyen2008a,Nguyen2008d} employs proper orthgogonal decomposition to generate the basis set and least-squares method to compute the interpolation point set. Generalized empirical interpolation method \cite{Maday2015a,Maday2013} generalizes the EIM concept by replacing the pointwise function evaluations by more general measures defined as linear functionals.

The first-order empirical interpolation method (FOEIM), introduced in \cite{Nguyen2023d}, enhances the original EIM by utilizing first-order partial derivatives of a parametrized nonlinear function to generate a larger set of basis functions and interpolation points. This method uses the same parameter sample set to construct these functions, thereby increasing the approximation power without requiring additional FOM solutions. FOEIM effectively improves the accuracy of hyper-reduced ROMs while maintaining computational efficiency. This paper extends the FOEIM introduced in \cite{Nguyen2023d} to time-dependent nonlinear PDEs. 


 \subsection{Empirical interpolation procedure}

We aim to interpolate the nonlinear function $g(u_N(\bm x, t_i, \bm \mu), \bm \mu)$ using a set of basis functions $\Psi_M^g = \mbox{span} \{\psi_m^g(\bm x), 1 \le m \le M\}$ and a set of interpolation points $T_M^g = \{{\bm x}_1^g, \ldots, {\bm x}_M^g\}$. The interpolant $g_M(\bm x, t_i, \bm \mu)$ is given by 
\begin{equation}
\label{eq1w}
g_M(\bm x, t_i, \bm \mu) = \sum_{m=1}^M \beta_{M,m}(t_i, \bm \mu) \psi^g_m(\bm x)   ,
\end{equation}
where the coefficients $\beta_{M,m}(t_i, \bm \mu), 1 \le m \le M,$ are found as the solution of the following linear system 
\begin{equation}
\label{eq2w}
\sum_{m=1}^M  \psi^g_m({\bm x}^g_p)   \beta_{M,m}(t_i,\bm \mu) = g(u_N({\bm x}_p^g, t_i, \bm \mu), \bm \mu), \quad 1 \le p \le M .
\end{equation}
It is convenient to compute the coefficient vector $\bm \beta_M(t_i,\bm \mu)$ as follows
\begin{equation}
\label{eqcoeff}
\bm \beta_M(t_i,\bm \mu) = \bm B^{-1}_M \bm b_M(t_i,\bm \mu),
\end{equation}
where $\bm B_M \in \mathbb{R}^{M \times M}$ has entries $B_{M,pm} =  \psi_m^g({\bm x}_p^g)$ and $\bm b_M(t_i,\bm \mu) \in \mathbb{R}^M$ has entries $b_{M,p}(t_i,\bm \mu) =  g(u_N({\bm x}_p^g, t_i, \bm \mu),  \bm \mu)$.

Next, we approximate ${g}_{N, j}(\bm \alpha_N(t_i,\bm \mu))$ defined in (\ref{eq:7-10}) with
\begin{equation}
\label{eq:7-10wq}
\widehat{g}_{N, j}(\bm \alpha_N(t_i,\bm \mu))  = \sum_{m=1}^M \beta_{M,m}(t_i, \bm \mu)  \int_{\Omega}   \psi^g_m(\bm x) \phi_j  .
\end{equation} 
Combining (\ref{eqcoeff}) and (\ref{eq:7-10wq}) yields 
\begin{equation}
\label{eq:7-10ww}
\widehat{\bm g}_{N}(\bm \alpha_N(t_i,\bm \mu))  = \bm G_{N M} \bm b_M(t_i,\bm \mu),
\end{equation} 
where $\bm G_{N M} = \bm H_{N M} \bm B_{M}^{-1}$ and $\bm H_{N M}$ has entries $H_{NM, jm } = \int_{\Omega}   \psi^g_m(\bm x) \phi_j$. The approximation accuracy depends crucially on both the subspace $\Psi^g_M$ and the interpolation point set $T^g_M$. In what follows, we discuss our approach to constructing both $\Psi^g_M$ and $T_M^g$.

\subsection{Basis functions}

In Section \ref{section2.1b}, we applied POD to the solution snapshot set $V^u_K = \{\zeta_k(\bm x) = u(\bm x, t_i, \bm \mu_j) \}_{k=1}^K$ to generate the RB space $\Phi_N$, in which the RB solution $u_N(t_i, \bm \mu)$ resides. Similarly, we can construct an RB space for the nonlinear function  $g(u_N(t_i, \bm \mu), \bm \mu)$. Assuming $g$ does not explicitly depend on $\bm \mu$, we define the snapshot set $V^g_K$ as follows
\begin{equation}
 V^g_K \equiv \{\xi_k(\bm x) = g(\zeta_k(\bm x)), 1 \le k \le K\} .   
\end{equation}
POD can then be applied to $V^g_K$ to generate the basis functions needed for approximating the nonlinear term $g(u_N(t_i, \bm \mu), \bm \mu)$. We generalize this method by including the partial derivatives of the nonlinear term.

The first-order empirical interpolation method generates basis functions by exploiting the first-order Taylor expansion of $g(w)$ at $v$:
\begin{equation}
G(w,v) = g(v) + g'_u(v) (w - v) ,
\end{equation}
where $g'_u$ denotes the partial derivative of $g(u)$ with respect to $u$. Taking $w = \zeta_k$ and $v = \zeta_{k'}$, where $(\zeta_k,\zeta_{k'})$ is any pair of two functions in $V_K^u$, we arrive at
\begin{equation}
\label{taylor}
G(\zeta_k,\zeta_{k'}) = g(\zeta_{k'}) + g'_u(\zeta_{k'}) (\zeta_{k} - \zeta_{k'})
\end{equation}
for $1 \le k,k' \le K$. We then define an expanded snapshot set
\begin{equation}
 V^g_{K^2} \equiv \{\rho_{kk'}(\bm x) = G(\zeta_k,\zeta_{k'}), 1 \le k,k' \le K\} .   
\end{equation}
Note that we have $V_K^g \subset V^g_{K^2}$. Furthermore, $V_{K^2}^{g}$ can be significantly  richer than $V_K^g$ as it contains considerably more snapshots. However, due to the potentially large number of snapshots in $V_{K^2}^{g}$, generating basis functions via POD becomes computationally intensive. Hence, we reduce the number of snapshots by using the nearest parameter points as follows.

First, we compute the distances $d_{jj'} = \|\bm \mu_j - \bm \mu_{j'}\|, 1 \le j,j' \le J,$ between each pair of parameter points $\bm \mu_{j}$  and  $\bm \mu_{j'}$ in $S_J$. For any given $\bm \mu_j \in S_J$, we define $S_L(\bm \mu_j)$ a set of exactly $L$ parameter points that are closest to $\bm \mu_j$. Then, for any snapshot $\zeta_k \in V_K^u$, we define $W_L(\zeta_k) = \{u(\bm x, t_i, \bm \mu_l), \, \forall \bm \mu_l \in S_L(\bm \mu_j) \}$, which contains  exactly $L$ snapshots. Next, we introduce the snapshot set
\begin{equation}
 V^g_{KL} \equiv \{\rho_{kk'}(\bm x) = G(\zeta_k,\zeta_{k'}), \ \forall \zeta_k \in V^u_K, \, \forall  \zeta_{k'} \in W_L(\zeta_k) \} .   
\end{equation}
This snapshot set is nested between $V_K^g$ and $V_{K^2}^g$, $V_K^g \subset V_{KL}^g \subset V^g_{K^2}$. For $L=1$, $V_{KL}^g$ reduces to $V_{K}^g$. For $L > 1$, the set $V_{KL}^g$ incorporates additional information from the partial derivatives. This inclusion of derivative information in 
$V_{KL}^g$ produces a more comprehensive snapshot set and improves the accuracy and stability of the resulting ROM.

Finally, we apply POD to the snapshot set $V^g_{KL}$ to generate $M$ basis functions, denoted as $\Phi_M^g =\{\phi_m^g, 1 \le m \le M \}$. The parameter $L$ only impacts the computational cost of POD during the offline stage and has no impact on other offline or online computations. This allows flexibility in choosing $L$ to balance accuracy and computational effort without affecting the efficiency of ROM in the online stage. However, the parameter $L$ can have a significant effect on the accuracy of the resulting ROM. We will study this effect by comparing the accuracy and stability of ROMs corresponding to different values of $L$.

\subsection{Interpolation points}

We multiply each basis function $\phi_m^g$ in $\Phi_M^g$ with $\sqrt{\lambda_m^g}$, where $\lambda^g_m$ is the corresponding eigenvalue. This is necessary for the following procedure to select the basis functions according to the magnitude of the eigenvalues. We apply the EIM procedure \cite{Maday2008c} to the space $\Phi_M^g$ to compute interpolation points and basis functions as follows.

The first interpolation point and basis function are given by
\begin{equation}
{\bm x}_1^g = \arg \sup_{\bm x \in \Omega} |\phi^g_{j_1}(\bm x)|, \qquad \psi^g_1(\bm x) = \phi^g_{j_1}(\bm x)/\phi^g_{j_1}({\bm x}^g_1),      
\end{equation}
where the index $j_1$ is determined by
\begin{equation}
\label{eq23w}
j_1 = \arg \max_{1 \le l \le M} \|\phi_l^g  \|_{L^\infty(\Omega)} .  
\end{equation}
For $m = 2, \ldots, M$, we solve the linear systems
\begin{equation}
\sum_{i=1}^{m-1}  \psi^g_i({\bm x}^g_j)   \sigma_{li} = \phi^g_l({\bm x}^g_j), \quad 1 \le j \le m-1 , 1 \le l \le M,
\end{equation}
we then find
\begin{equation}
\label{argmax}
j_{m} = \arg \max_{1 \le l \le M} \| \phi^g_l(\bm x) - \sum_{i=1}^{m-1}  \sigma_{li} \psi^g_i(\bm x)\|_{L^\infty(\Omega)}, 
\end{equation}
and set
\begin{equation}
{\bm x}^g_m = \arg \sup_{\bm x \in \Omega} |r_M(\bm x)|, \qquad \psi^g_m(\bm x) = r_m(\bm x)/r_m({\bm x}^g_m)    ,   
\end{equation}
where the residual function $r_m(\bm x)$ is given by
\begin{equation}
\label{eq28w}
r_m(\bm x) = \phi^g_{j_m}(\bm x) - \sum_{i=1}^{m-1}  \sigma_{j_m i} \psi^g_i(\bm x) . 
\end{equation}
In practice, the supremum $\sup_{\bm x \in \Omega} |r_M(\bm x)|$ is computed on the set of quadrature points on all elements in the mesh. In other words, the interpolation points $\{\bm x^g_m\}_{m=1}^M$ are selected from the quadrature points.

\subsection{Error estimate of the first-order empirical interpolation method}

The first-order EIM constructs the interpolation point set $T_M^g = \{\bm x_m^g\}_{m=1}^M$ and the basis set $\Psi^g_M = \mbox{span}\{\psi_m^g\}_{m=1}^M$ by leveraging both the solution snapshots in $V_N^u$ and their partial derivatives. The first-order EIM is well-defined in the sense that the basis functions are linearly independent and that the interpolation procedure yields a unique interpolant. Furthermore,  the following error estimate  provides a measure of the interpolation error.   
\begin{thm}
\label{thm3}
If $g(u_N(\bm x, t, \bm \mu)) \in \Psi_{M+P}$ for $P \in \mathbb{N}_{+}$, then the interpolation error $\varepsilon_M(t,\bm \mu) \equiv \|g(u_N(\bm x, t, \bm \mu)) -  g_M(\bm x, t, \bm \mu) \|_{L^\infty(\Omega)}$ is bounded by
\begin{equation}
\label{EIMbound3}
\varepsilon_M(t,\bm \mu) \le \hat{\varepsilon}_{M,P} (t,\bm \mu) \equiv  \sum_{j=1}^P | e_j(t,\bm \mu) | ,
\end{equation}
where $e_j(t,\bm \mu), 1 \le j \le P,$ solve the following linear system
\begin{equation}
\sum_{j=1}^{P} \psi^g_{M+j}({\bm x^g}_{M+i}) e_j(t,\bm \mu)  =   g (u({\bm x}^g_{M+i},t, \bm \mu)) - g_M 
({\bm x}^g_{M+i}, t, \bm \mu ), \quad 1 \le i \le P .
\end{equation}
\end{thm}

The proof of this theorem can be found in \cite{Nguyen2023d,Nguyen2024}. The operation count of evaluating the error estimator (\ref{EIMbound3}) is only $O(P^2)$. Hence, the error estimator is very inexpensive to evaluate.  The error estimate  provides a form of certification, allowing us to control and balance the trade-off between computational efficiency and accuracy.

\subsection{A simple test case}

\begin{figure}[htbp]
	\centering
	\begin{subfigure}[b]{0.49\textwidth}
		\centering		\includegraphics[width=\textwidth]{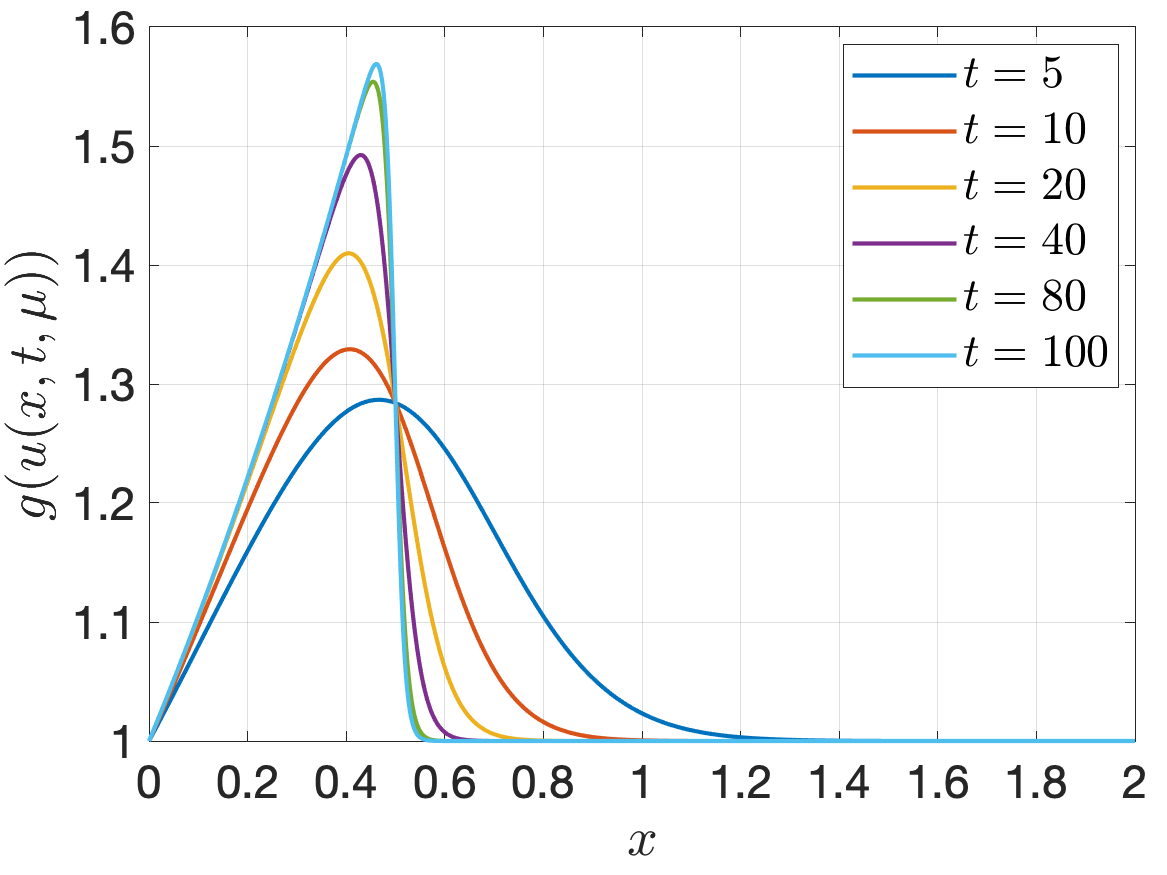}
		\caption{$\mu = 0$.}
	\end{subfigure}
	\hfill
	\begin{subfigure}[b]{0.49\textwidth}
		\centering		\includegraphics[width=\textwidth]{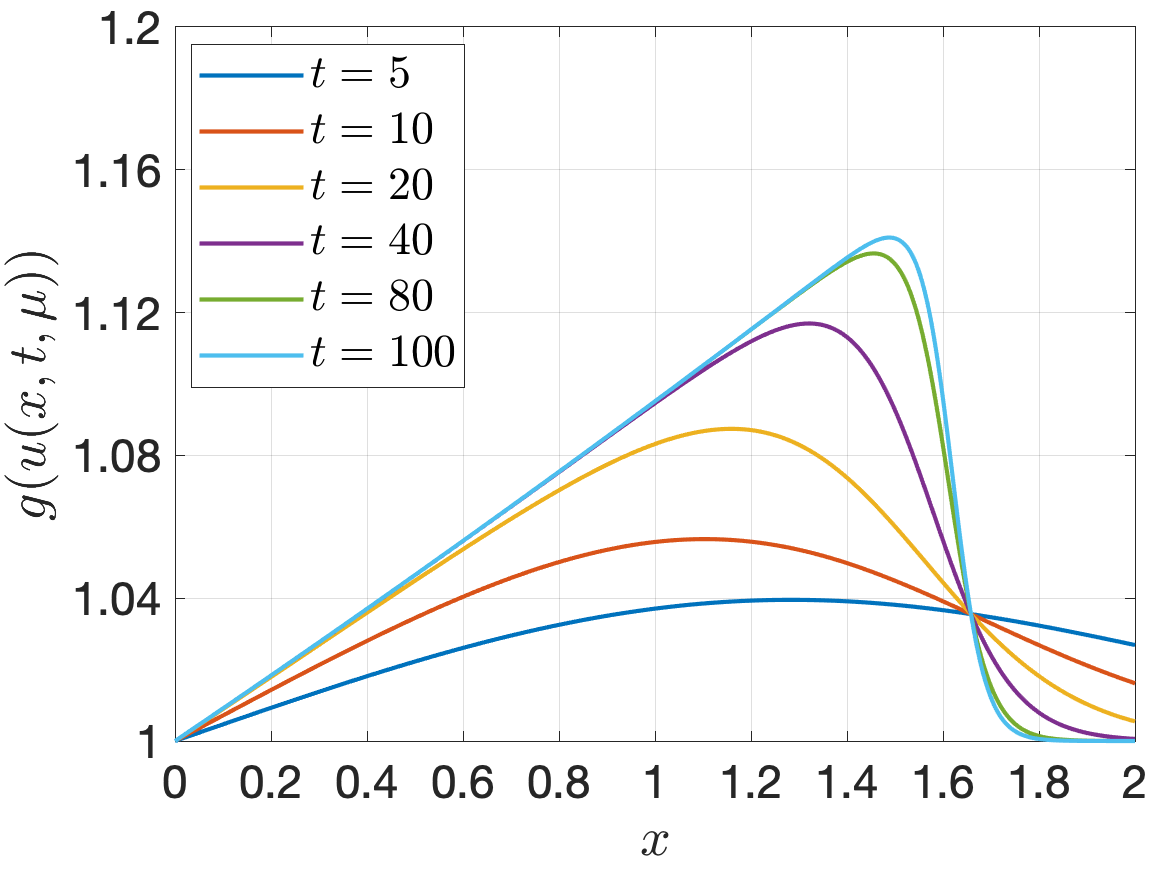}
		\caption{$\mu = 10$.}
	\end{subfigure}
	\caption{Plots of $g(u(x,t,\mu))$ as a function of $x$ and $t$ for $\mu = 0$ and $\mu =10$.}
	\label{ex1fig1}
\end{figure}

We present numerical results from a simple test case to study the performance of the first-order EIM. The test case involves the following parametrized functions 
$$u(x, t, \mu) = \frac{x}{(\mu + 1) \left(1 + \sqrt{\frac{\mu+1}{\exp(0.5 t)}} \exp \left(\frac{t x^2}{\mu+1} \right) \right) }, \qquad g(u) = \exp(u)$$
in a physical domain $\Omega = [0,2]$, time domain $T = [0, 100]$, and parameter domain $\mathcal{D} = [0, 10]$. Figure \ref{ex1fig1} shows instances of the nonlinear function $g$ for two different values of $\mu$. 


\begin{figure}[htbp]
	\centering
	\begin{subfigure}[b]{0.49\textwidth}
		\centering		\includegraphics[width=\textwidth]{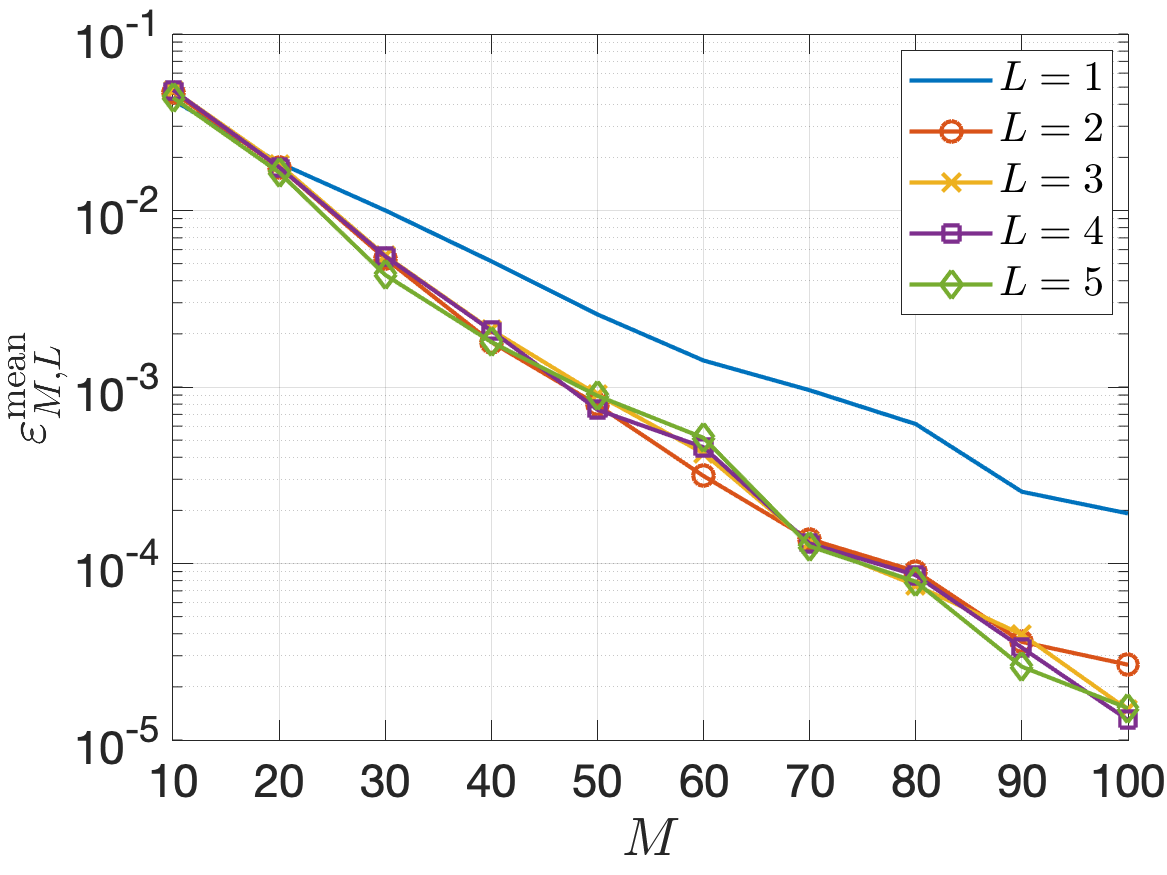}
		\caption{$J = 6$.}
	\end{subfigure}
	\hfill
	\begin{subfigure}[b]{0.49\textwidth}
		\centering		\includegraphics[width=\textwidth]{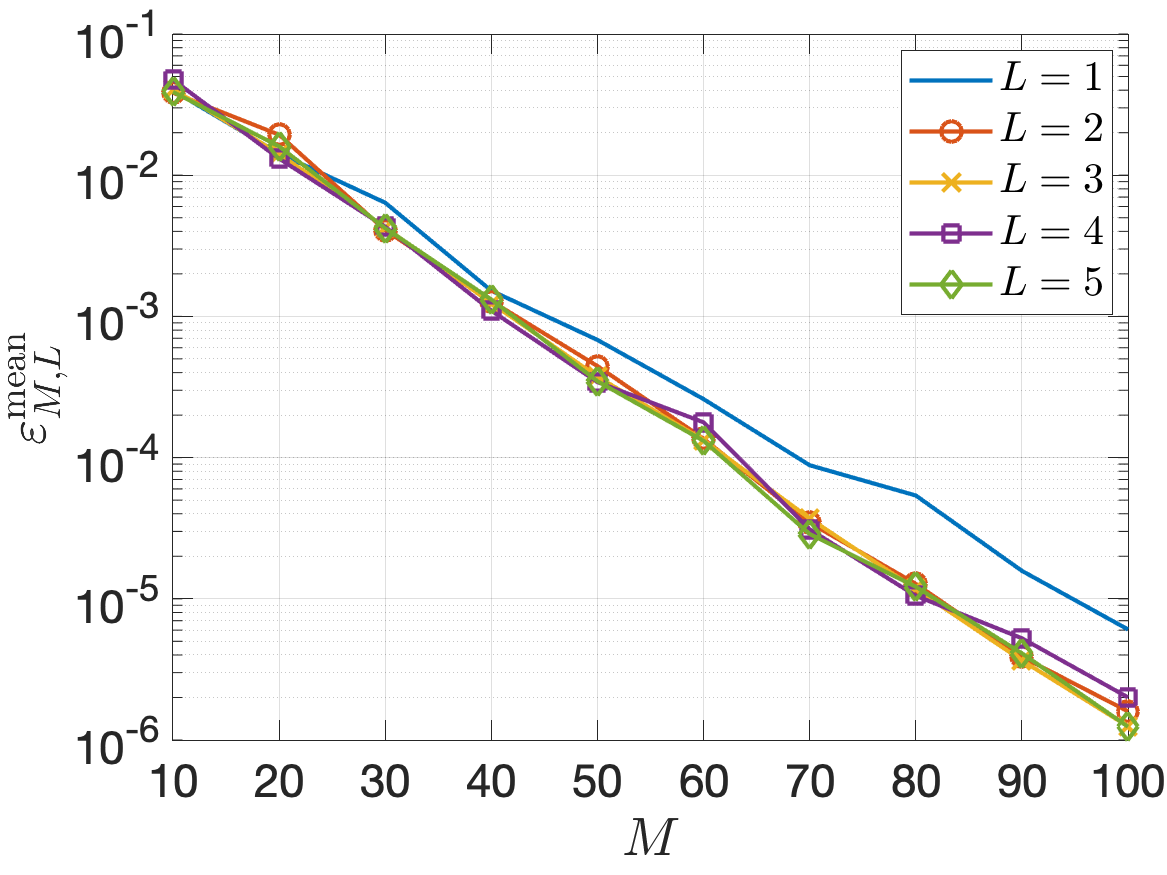}
		\caption{$J = 12$.}
	\end{subfigure}
	\caption{Mean absolute error as a function of $M$ and $L$ for $J = 6$ and $J =12$.}
	\label{ex1fig2}
\end{figure}

We consider $S_{J_{\max}} = \{0, 10, 1.4, 8.6, 4.2, 5.8, 0.5, 9.5, 2.7, 7.3, 0.15, 9.85 \}$ to sample the parameter domain, where $J_{\max} = 12$ is the  maximum number of parameter points. The parameter samples are hierarchical such that $S_J \subset S_{J_{\max}}$ for any $J \in [1, J_{\max}]$. We discretize the time domain into $I=100$ equal intervals and the physical domain into 1000 uniform elements. We use a test sample $\Xi^{\rm test}$ of $N^{\rm test} = 100$ parameter points distributed uniformly in the parameter domain. The mean error is defined as $\varepsilon_{M,L}^{\rm mean} = \frac{1}{I \, N^{\rm test}}\sum_{\mu \in \Xi^{\rm test}} \sum_{i=1}^I \varepsilon_M(t_i,\mu)$, where $\varepsilon_M(t_i,\mu)$ represents the interpolation error. Figure \ref{ex1fig2} plots $\varepsilon_{M,L}^{\rm mean}$ as a function of both $M$ and $L$ for $J = 6$ and $J = 12$. These plots illustrate how FOEIM reduces the average error as both $M$ and $L$ increase. The error decreases as $M$ increases, indicating that adding more basis functions improves the approximation accuracy. When $L = 1$, the error decreases more slowly compared to higher values of $L$, showing that increasing $L$ enhances the convergence rate of the error. For $L \geq 2$, the curves almost overlap, suggesting diminishing returns in accuracy improvement as $L$ increases beyond 2. Hence, incorporating the partial derivatives of the nonlinear function significantly improves interpolation accuracy, particularly when the number of parameter points is limited. 

We define the mean error estimate as $\widehat{\varepsilon}_{M,L}^{\rm mean} = \frac{1}{I \, N^{\rm test}}\sum_{\mu \in \Xi^{\rm test}} \sum_{i=1}^I \widehat{\varepsilon}_{M,P}(t_i,\mu)$ and the mean effectivity as $\widehat{\eta}_{M,L}^{\rm mean} = \frac{1}{I\, N^{\rm test}}\sum_{\mu \in \Xi^{\rm test}} \sum_{i=1}^I \widehat{\varepsilon}_{M,P}(t_i,\mu)/\varepsilon_M(t_i,\mu)$. Tables \ref{ex1tab2} and \ref{ex1tab3} present the values of $\widehat{\varepsilon}_{M,P}^{\rm mean}$ and $\widehat{\eta}_{M,P}^{\rm mean}$ as a function $M$ and $L$ for $J=6$ and $J=12$, respectively. 
The results indicate that the error estimates decrease rapidly as $M$ increase. Additionally, the error estimates are highly accurate, as reflected by the mean effectivities being consistently less than 3. This demonstrates that the estimation method provides sharp bounds on the actual error, ensuring efficiency, accuracy and reliability in the interpolation.

\begin{table}[htbp]
\centering
\small
	\begin{tabular}{|c||cc|cc|cc|cc|}
		\cline{1-7}
    &		 
	 \multicolumn{2}{|c|}{$L=1$} & \multicolumn{2}{c|}{$L=2$} & 
		 \multicolumn{2}{c|}{$L=3$} \\   
   $M$ & $\widehat{\varepsilon}_{M,L}^{\rm mean}$ & $\widehat{\eta}_{M,L}^{\rm mean}$ & $\widehat{\varepsilon}_{M,L}^{\rm mean}$ & $\widehat{\eta}_{M,L}^{\rm mean}$ & $\widehat{\varepsilon}_{M,L}^{\rm mean}$ & $\widehat{\eta}_{M,L}^{\rm mean}$ \\
		\cline{1-7}
  10  &  6.26\mbox{e-}2  &  1.52  &  6.57\mbox{e-}2  &  1.39  &  8.00\mbox{e-}2  &  1.29  \\  
  20  &  2.03\mbox{e-}2  &  1.45  &  2.69\mbox{e-}2  &  1.60  &  3.53\mbox{e-}2  &  1.75  \\  
  30  &  1.64\mbox{e-}2  &  1.53  &  8.10\mbox{e-}3  &  1.80  &  7.29\mbox{e-}3  &  1.63  \\  
  40  &  8.49\mbox{e-}3  &  1.38  &  2.92\mbox{e-}3  &  1.67  &  3.96\mbox{e-}3  &  1.99  \\  
  50  &  3.87\mbox{e-}3  &  1.93  &  9.51\mbox{e-}4  &  1.80  &  1.29\mbox{e-}3  &  2.09  \\  
  60  &  1.01\mbox{e-}3  &  1.83  &  5.44\mbox{e-}4  &  1.96  &  6.04\mbox{e-}4  &  1.77  \\  
  70  &  2.30\mbox{e-}3  &  1.93  &  2.84\mbox{e-}4  &  2.19  &  1.86\mbox{e-}4  &  1.90  \\  
  80  &  9.64\mbox{e-}4  &  2.29  &  1.40\mbox{e-}4  &  2.14  &  1.85\mbox{e-}4  &  2.26  \\  
  90  &  6.60\mbox{e-}4  &  2.02  &  5.58\mbox{e-}5  &  2.29  &  8.84\mbox{e-}5  &  2.98  \\  
  100  &  3.46\mbox{e-}4  &  1.93  &  7.24\mbox{e-}5  &  2.33  &  2.29\mbox{e-}5  &  2.38  \\  
\hline
		
	\end{tabular}
	\caption{ $\widehat{\varepsilon}_{M,L}^{\rm mean}$ and $\widehat{\eta}_{M,L}^{\rm mean}$ as a function of $M$ and $L$ for $P = 5$ and $J=6$.} 
	\label{ex1tab2}
\end{table}

\begin{table}[htbp]
\centering
\small
	\begin{tabular}{|c||cc|cc|cc|cc|}
		\cline{1-7}
    &		 
	 \multicolumn{2}{|c|}{$L=1$} & \multicolumn{2}{c|}{$L=2$} & 
		 \multicolumn{2}{c|}{$L=3$} \\   
   $M$ & $\widehat{\varepsilon}_{M,L}^{\rm mean}$ & $\widehat{\eta}_{M,L}^{\rm mean}$ & $\widehat{\varepsilon}_{M,L}^{\rm mean}$ & $\widehat{\eta}_{M,L}^{\rm mean}$ & $\widehat{\varepsilon}_{M,L}^{\rm mean}$ & $\widehat{\eta}_{M,L}^{\rm mean}$ \\
		\cline{1-7}
    10  &  3.71\mbox{e-}2  &  1.29  &  3.52\mbox{e-}2  &  1.14  &  3.40\mbox{e-}2  &  1.09  \\  
  20  &  1.77\mbox{e-}2  &  1.35  &  3.71\mbox{e-}2  &  1.64  &  2.01\mbox{e-}2  &  1.01  \\  
  30  &  9.65\mbox{e-}3  &  1.70  &  4.65\mbox{e-}3  &  1.31  &  5.37\mbox{e-}3  &  1.22  \\  
  40  &  2.00\mbox{e-}3  &  1.52  &  2.11\mbox{e-}3  &  1.89  &  2.68\mbox{e-}3  &  2.14  \\  
  50  &  1.04\mbox{e-}3  &  1.62  &  6.98\mbox{e-}4  &  1.79  &  5.85\mbox{e-}4  &  1.68  \\  
  60  &  4.32\mbox{e-}4  &  1.71  &  2.42\mbox{e-}4  &  1.84  &  2.44\mbox{e-}4  &  1.72  \\  
  70  &  1.45\mbox{e-}4  &  1.83  &  5.81\mbox{e-}5  &  1.65  &  4.42\mbox{e-}5  &  1.89  \\  
  80  &  1.18\mbox{e-}4  &  2.10  &  1.98\mbox{e-}5  &  2.01  &  2.61\mbox{e-}5  &  2.14  \\  
  90  &  2.70\mbox{e-}5  &  1.77  &  5.00\mbox{e-}6  &  2.05  &  6.80\mbox{e-}6  &  1.87  \\  
  100  &  7.97\mbox{e-}6  &  2.40  &  3.66\mbox{e-}6  &  2.39  &  2.53\mbox{e-}6  &  2.62  \\  
\hline
	\end{tabular}
	\caption{$\widehat{\varepsilon}_{M,L}^{\rm mean}$ and $\widehat{\eta}_{M,L}^{\rm mean}$ as a function of $M$ and $L$ for $P=5$ and $J=12$.} 
	\label{ex1tab3}
\end{table}

\section{Numerical Experiments}
\label{section3.4}

In this section, we present numerical results from two parametrized nonlinear PDEs to demonstrate the hyperreduced Galerkin-Newton method, specifically the FOEIM-GN (First-Order Empirical Interpolation Method--Galerkin-Newton) approach. The results focus on comparing the performance of the FOEIM--GN method against the standard Galerkin-Newton (GN) method and the FOM. These comparisons will highlight the accuracy, computational efficiency, and potential advantages of using the FOEIM-GN method in reducing the complexity of solving time-dependent nonlinear PDEs. The PDEs studied include the Allen-Cahn equation, which is central to phase separation in multi-component alloy systems, and the Buckley-Leverett equation, which models two-phase flow in porous media. The Allen-Cahn equation is particularly important for understanding phase transitions in materials science and curvature-driven flows in geometry, while the Buckley-Leverett equation is crucial for predicting fluid displacement fronts in immiscible displacement processes.

Both GN and FOEIM-GN are evaluated on the basis of accuracy, computational cost, and convergence behavior. To assess the accuracy, we define the following  errors  
\begin{equation}
\epsilon^u_{N}(t_i,\bm \mu) = \|u(t_i,\bm \mu) - u_{N}(t_i,\bm \mu)\|_X, \quad \epsilon^s_{N}(t_i,\bm \mu) = |s(t_i,\bm \mu) - s_{N}(t_i,\bm \mu)| ,
\end{equation}
and the average errors
\begin{equation}
\bar{\epsilon}_{N}^u =  \frac{1}{I N^{\rm test}}\sum_{\bm \mu \in \Xi^{\rm test}} \sum_{i=1}^I \epsilon_N^u (t_i,\bm \mu), \quad  \bar{\epsilon}_{N}^s = \frac{1}{I N^{\rm test}} \sum_{\bm \mu \in \Xi^{\rm test}} \sum_{i=1}^I \epsilon_N^s (t_i,\bm \mu)  .  
\end{equation}  
where $\Xi^{\rm test}$ is a test sample of $N^{\rm test}$ parameter points  distributed uniformly in the parameter domain.


\subsection{Buckley-Leverett equation}

We consider the two-dimensional Buckley-Leverett equation consisting of  linear diffusion and nonlinear convection terms
\begin{equation}
\dot{u}^{\rm e} -\mu \nabla^2u^{\rm e} + \nabla \cdot \bm f(u^{\rm e}) = 0, \quad \mbox{in } \Omega \times (0, T), 
\end{equation} 
 with initial condition $u^{\rm e}(\bm x, 0, \mu) = \exp(-16 (x_1^2 + x_2^2))$ and homogeneous Dirichlet condition on the boundary $\partial \Omega$, where $\Omega = (-1.5,1.5)^2$, $T = 1$, and $\mu  \in {\cal D} \equiv [0.03, 0.1]$. The flux vector $\bm f(u) = (f^1(u), f^2(u))$ consists of non-linear functions of the field variable $u$ as $f^1(u) = \frac{u^2}{u^2 + (1-u)^2}$ and $f^2(u) = f^1(u)(1 - 5 (1-u)^2)$. The output of interest is the average of the field variable over the
physical domain. 

The weak formulation is then stated as: given $\mu \in
{\cal D}$, we find $s(t,\mu) = \int_\Omega u(t,\mu)$, where $u(t, \mu) \in X \subset H_0^1(\Omega) \equiv \{v \in
H^1(\Omega) \mbox{ } | \mbox{ } v|_{\partial \Omega} = 0\}$ is the solution
of
\begin{equation}
m(\dot{u}(t, \mu), v) + \mu a^0(u(t, \mu), v) + b (u(t, \mu), v) = 0, \quad \forall  v
\in X \ , 
\end{equation}
where
\begin{equation}
m(w, v) = \int_\Omega  w v, \quad a^0(w, v) = \int_\Omega \nabla w \cdot \nabla v, \quad  b(w,v) = - \int_\Omega \bm f(w) \cdot \nabla v .
\end{equation}
The finite element (FE) approximation space is $X = \{v \in H_0^1(\Omega) : v|_K \in \mathcal{P}^3(T), \  \forall T \in \mathcal{T}_h \}$, where $\mathcal{P}^3(T)$ is a space of polynomials of degree $2$ on an element $T \in \mathcal{T}_h$ and $\mathcal{T}_h$ is a finite element grid of $80 \times 80$ quadrilaterals. The dimension of the FE space is $\mathcal{N} = 9409$. For the temporal discretization, we use the Backward-Euler scheme with timestep size $\Delta t = 0.01$. The training sample $S_{J} = \{0.03,    0.044,    0.058,    0.072,    0.086,    0.1\}$ has 6 parameter points, while the test sample $\Xi_{\rm test}$ consists of $11$ parameter points distributed uniformly in the parameter domain.

\begin{figure}[htbp]
	\centering
	\begin{subfigure}[b]{0.49\textwidth}
		\centering		\includegraphics[width=\textwidth]{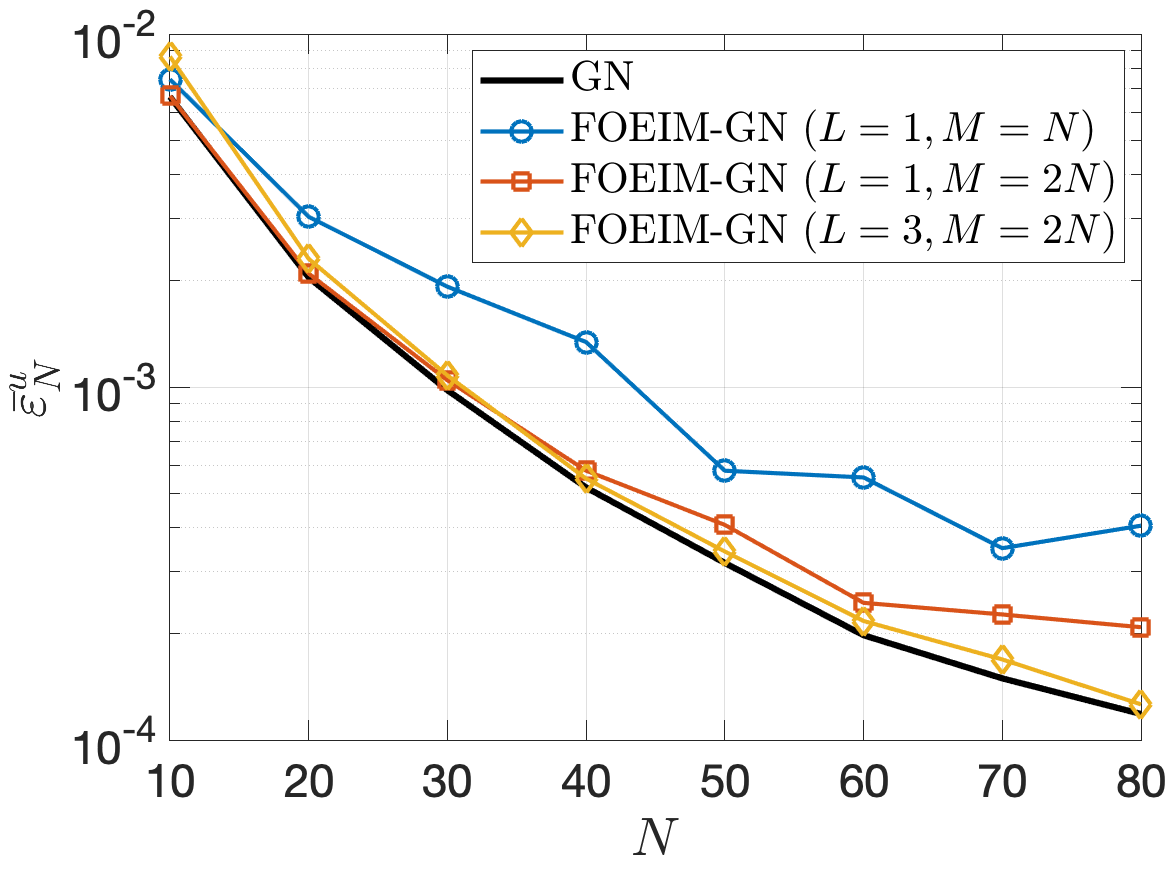}
		\caption{Average solution error $\bar{\epsilon}_{N}^u$.}
	\end{subfigure}
	\hfill
	\begin{subfigure}[b]{0.49\textwidth}
		\centering		\includegraphics[width=\textwidth]{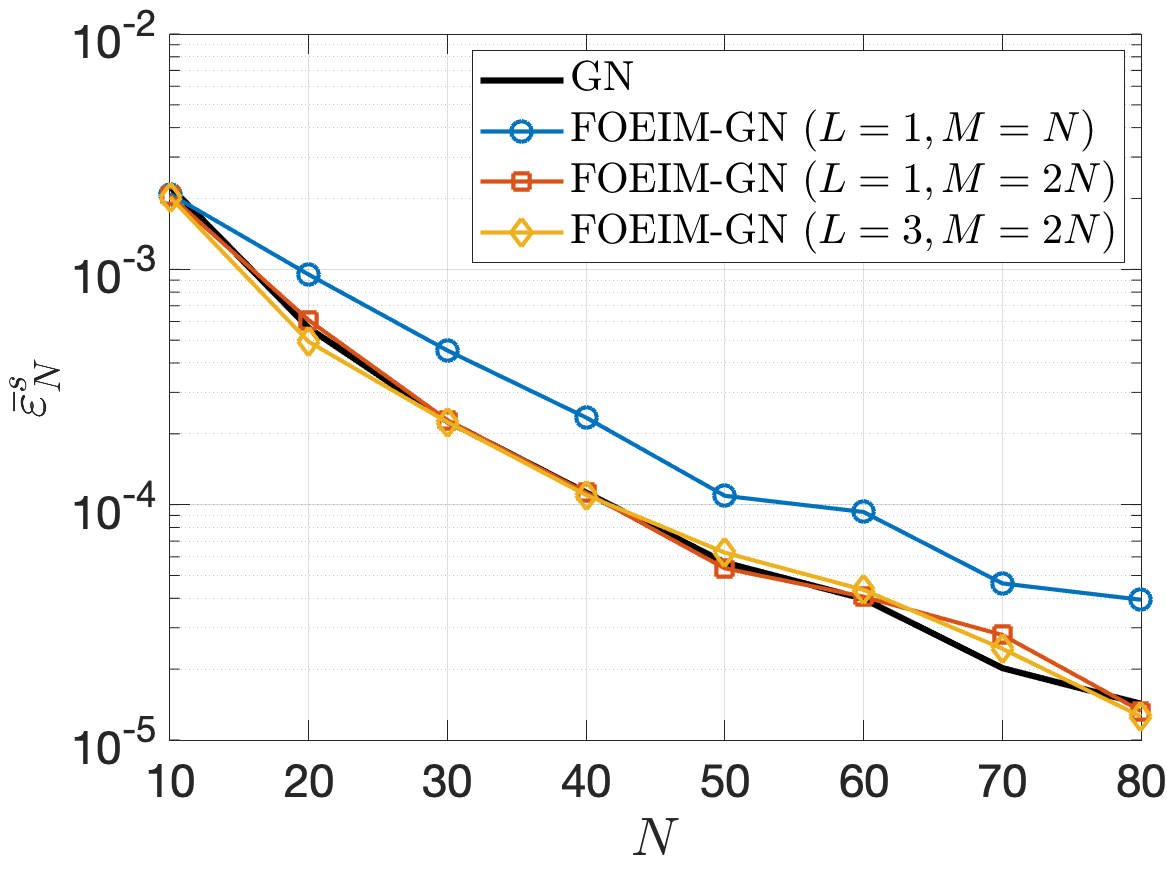}
		\caption{Average output error $\bar{\epsilon}_{N}^s$.}
	\end{subfigure}
	\caption{Comparison of accuracy between the GN method and the FOEIM-GN method with three different values of $L$ and $M$.}
	\label{ex2fig2}
\end{figure}

Figure \ref{ex2fig2} presents the convergence of the mean solution error $\bar{\epsilon}_{N}^u$  and the mean output error $\bar{\epsilon}_{N}^s$ as functions of $N$ for the GN method and the FOEIM-GN method with three different values of $L$ and $M$. As $N$ increases, the methods exhibit a clear reduction in error. Increasing $M$ from $N$ to $2N$ leads to a significant reduction in error. When $L = 1$, the error decreases more slowly compared to $L = 3$, showing that increasing $L$ improves the convergence rate of the error. Hence, incorporating the partial derivatives of the nonlinear terms enhances the accuracy of the FOEIM-GN method. For $L=3$ and $M=2N$, the FOEIM-GN method closely mirrors the accuracy of the standard GN method.

Table \ref{ex2tab2} shows the computational speedup for the GN method and the FOEIM-GN method compared to the finite element method (FEM) for different values of $N$.  Across all values of $N$, the GN method achieves a modest speedup (around 1.5–2.0x) compared to FEM. In contrast, the FOEIM-GN method is two and three orders of magnitude faster than FEM, although its speedup factor decreases as $N$ increases. For $M=N$, its speedup factor ranges from 5000x at $N=10$ to about 340x at $N=80$. For $M=2N$, its speedup factor ranges from 3310x at $N=10$ to about 280x at $N=80$. Thus, the FOEIM-GN method is two or three orders of magnitude faster than the GN method while achieving a similar accuracy.


\begin{table}[htbp]
\centering
\small
	\begin{tabular}{|c||c|c|c|}
		\cline{1-4}
  $N$  & \mbox{ } GN \mbox{ } & FOEIM-GN ($M=N$)& FOEIM-GN ($M=2N$) \\     
		\cline{1-4}
 10  &  1.89  &  5001.52  &  3309.12  \\  
 20  &  1.69  &  2022.39  &  1276.66  \\  
 30  &  1.67  &  984.65  &  817.89  \\  
 40  &  1.59  &  728.84  &  661.78  \\  
 50  &  1.52  &  574.33  &  599.62  \\  
 60  &  1.51  &  477.62  &  410.27  \\  
 70  &  1.49  &  390.45  &  333.23  \\  
 80  &  1.47  &  339.11  &  283.24  \\   
		\hline
	\end{tabular}
	\caption{Computational speedup  relative to the finite element method (FEM)  for the GN and FOEIM-GN methods as a function of $N$. The speedup is calculated as the ratio between the computational time of FEM and the online computational time of ROM.} 
	\label{ex2tab2}
\end{table}

\subsection{Allen-Cahn equation}

The Allen–Cahn (AC) equation is a reaction-diffusion equation composed of the reaction term and the diffusion term
\begin{equation}
\dot{u}^{\rm e} - \nabla^2u^{\rm e} + \frac{g(u^{\rm e}) - u^{\rm e}}{\varepsilon^2} = 0, \quad \mbox{in } \Omega \times (0, T), 
\end{equation} 
with homogeneous Neumann condition $\nabla u^{\rm e} \cdot \bm n = 0$ on the boundary $\partial \Omega$, where $\Omega = (0,1)^2$ and $T = 0.02$.  Here the nonlinear term $g(u^{\rm e}) = (u^{\rm e})^3$ originates from the derivative of a potential energy function, and $\varepsilon = 0.015$ is the thickness of the transition layer which is a small positive constant value. The quantity $u^{\rm e}$ is an order parameter, which is one of the concentrations of the two components in a binary mixture. For example, $u^{\rm e} = 1$ in the one phase and $u^{\rm e} = 0$ in the other phase.  The initial condition is a parameter-dependent star shape 
\begin{equation}
u^{\rm e}(\bm x, 0, \mu) = \tanh \left( \frac{\mu + 0.1 \cos(6 \theta) - \sqrt{(x_1-0.5)^2 + (x_2-0.5)^2}}{\sqrt{2} \varepsilon} \right)
\end{equation}
with $\theta = \tan^{-1}\left(\frac{x_2-0.5}{x_1-0.5}\right)$ and $\mu  \in {\cal D} \equiv [0.25, 0.35]$. The parameter $\mu$ determines the area of the star shape. The output of interest is the average of the field variable over the physical domain.

The weak formulation is then stated as: given $\mu \in
{\cal D}$, we find $s(t,\mu) = \int_\Omega u(t,\mu)$, where $u(t, \mu) \in X \subset H^1(\Omega)$ is the solution
of
\begin{equation}
m(\dot{u}(t, \mu), v) +  a^0(u(t, \mu), v) + b (u(t, \mu), v) = 0, \quad \forall  v
\in X \ , 
\label{eq:7-6}
\end{equation}
where
\begin{equation}
m(w, v) = \int_\Omega  w v, \quad a^0(w, v) = \int_\Omega \nabla w \cdot \nabla v - \frac{1}{\varepsilon}\int_\Omega  w v, \quad  b(w,v) =  \frac{1}{\varepsilon} \int_\Omega g(w) v .
\label{eq:7-6a}
\end{equation}
The finite element (FE) approximation space is $X = \{v \in H^1(\Omega) : v|_K \in \mathcal{P}^3(T), \  \forall T \in \mathcal{T}_h \}$, where $\mathcal{P}^3(T)$ is a space of polynomials of degree $2$ on an element $T \in \mathcal{T}_h$ and $\mathcal{T}_h$ is a finite element grid of $80 \times 80$ quadrilaterals. The dimension of the FE space is $\mathcal{N} = 9409$. For the temporal discretization, we use the Backward-Euler scheme with timestep size $\Delta t = 0.0001$. The training sample $S_{J} = \{0.25,   0.27,  0.29,  0.31,  0.33,   0.35\}$ has 6 parameter points, while the test sample $\Xi_{\rm test}$ consists of $11$ parameter points distributed uniformly in the parameter domain.

\begin{figure}[htbp]
	\centering
	\begin{subfigure}[b]{0.49\textwidth}
		\centering		\includegraphics[width=\textwidth]{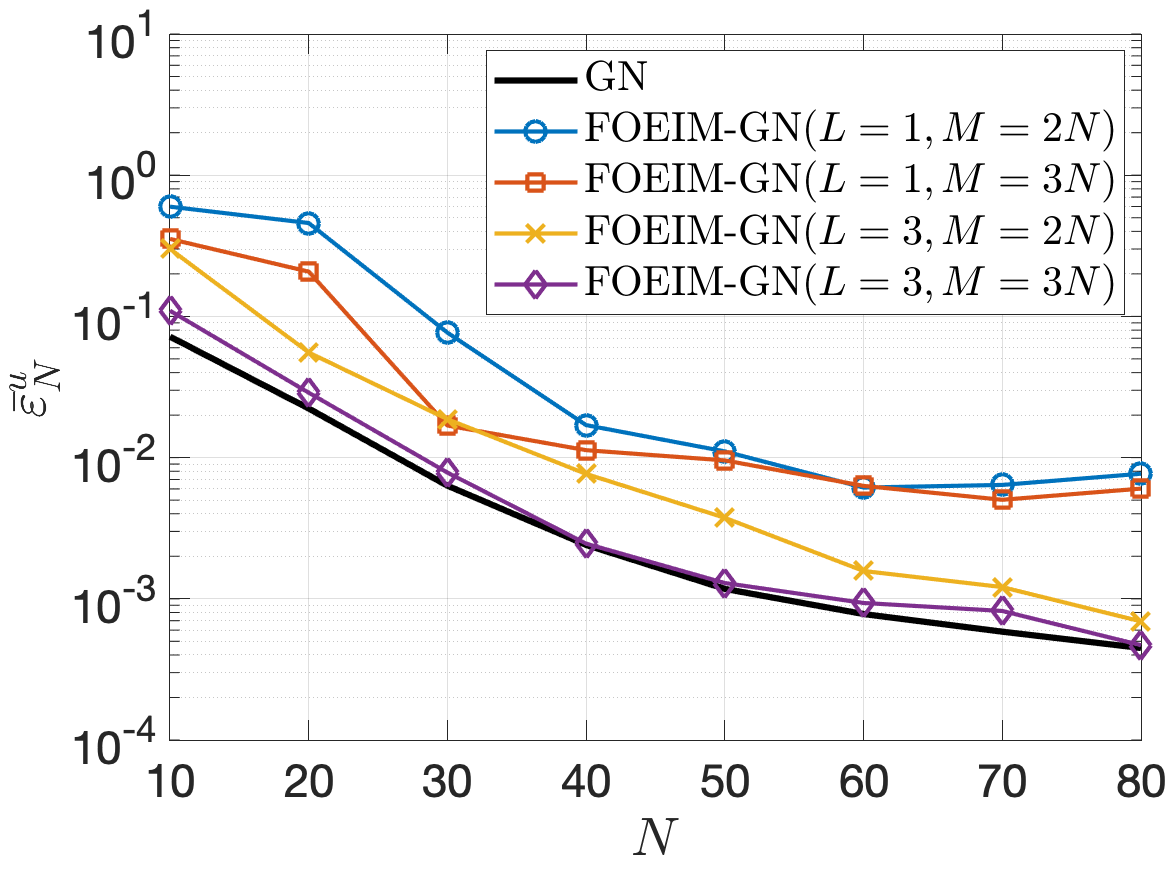}
		\caption{Average solution error $\bar{\epsilon}_{N}^u$.}
	\end{subfigure}
	\hfill
	\begin{subfigure}[b]{0.49\textwidth}
		\centering		\includegraphics[width=\textwidth]{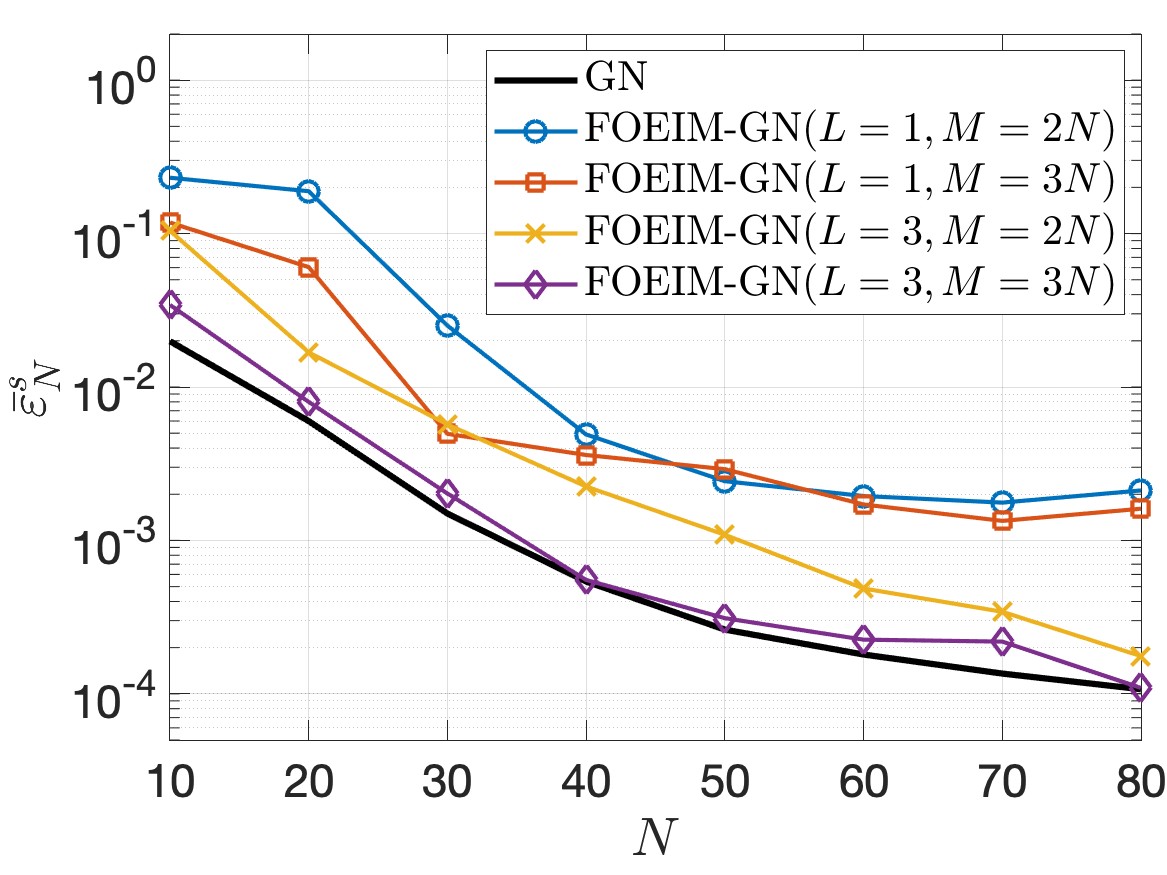}
		\caption{Average output error $\bar{\epsilon}_{N}^s$.}
	\end{subfigure}
	\caption{Comparison of accuracy between the GN method and the FOEIM-GN method with four different values of $L$ and $M$.}
	\label{ex3fig2}
\end{figure}

Figure \ref{ex3fig2} presents the convergence of the mean solution error $\bar{\epsilon}_{N}^u$  and the mean output error $\bar{\epsilon}_{N}^s$ as functions of $N$ for the GN method and the FOEIM-GN method with four different values of $L$ and $M$. The GN method consistently shows a steady decrease in error as 
$N$ increases. The FOEIM-GN methods exhibit varying degrees of accuracy improvement depending on the values of $L$ and $M$. For the FOEIM-GN method, $L=1$ and $M=2N$ yields higher errors than the other values of $L$ and $M$. Increasing $M$ from $2N$ to $3N$  results in smaller errors. Increasing  $L$ from $1$ to $3$ considerably improves accuracy. The  FOEIM-GN method with $L=3$ and $M=3N$ yields the smallest errors among all the four cases. For $L=3$ and $M=3N$, the FOEIM-GN method has almost the same errors as the standard GN method.

\begin{figure}[h!]
	\centering
\includegraphics[width=0.99\textwidth]{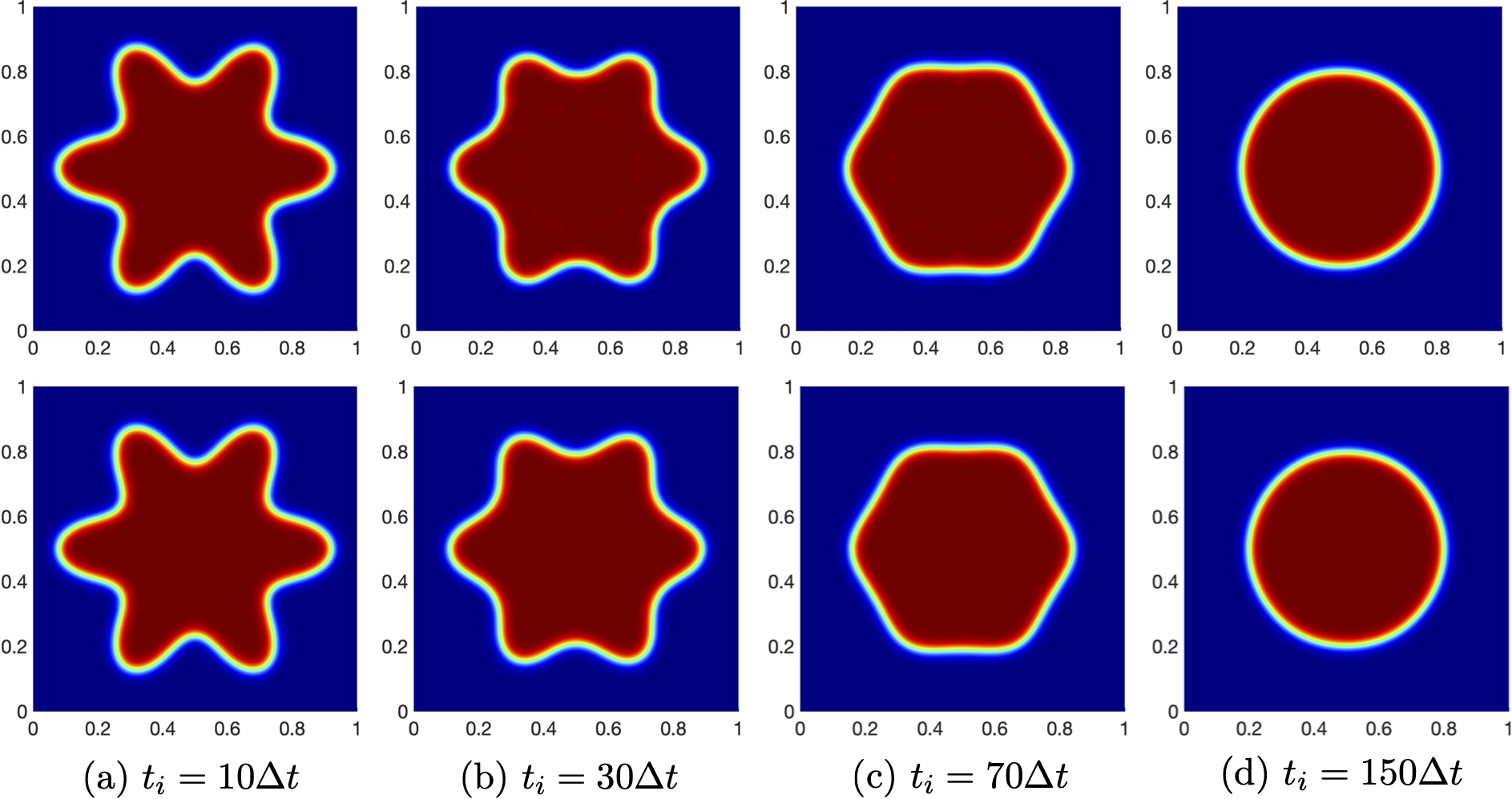}
 \caption{Evolution of the FOEIM-GN solution  with $N = 50, M = 100$, and $L=1$ (top row) or $L = 3$ (bottom row) for $\mu = 0.34$.}
	\label{ex3fig3}
\end{figure}

Figure \ref{ex3fig3} illustrates the evolution of the  FOEIM-GN solutions  with $N = 50, M = 100$, and $L=1$ (top row) or $L = 3$ (bottom row) for $\mu = 0.34$ at different time steps.  Initially, both configurations capture a star-shaped pattern of the initial solution, representing the early phase separation, with little difference between them. As time progresses to intermediate steps, the tips of the star move inward and the gaps between the tips move outward. The pattern begins to smoothen, characteristic of the Allen-Cahn dynamics, with the $L=3$ configuration showing a slightly smoother transition in the interface. By the time step $t_i=150 \Delta t$, both configurations evolve into a nearly perfect circle, representing the system's approach to equilibrium, though the $L=3$ configuration provides a more symmetric and smoother representation. Once the form deforms to a circular shape, the radius of the circle shrinks with increasing speed. The results suggest that increasing $L$ improves the ability to accurately capture the dynamics of phase separation, particularly at later stages of evolution.

Table \ref{ex3tab2} shows the computational speedup for the GN and FOEIM-GN methods compared to FEM for different values of $N$. GN is only about 2 times faster than FEM. In contrast,  FOEIM-GN is two and three orders of magnitude faster than FEM, although its speedup factor decreases as $N$ increases. For FOEIM-GN, the speedup factor ranges from 4000x at $N=10$ to about 300x at $N=80$. Thus,  FOEIM-GN is two or three orders of magnitude faster than GN for the same level of accuracy.

\begin{table}[h!]
\centering
\small
	\begin{tabular}{|c||c|c|c|}
		\cline{1-4}
  $N$  & \mbox{ } GN \mbox{ } & FOEIM-GN ($M=2N$)& FOEIM-GN ($M=3N$) \\     
		\cline{1-4}
 10  &  2.29  &  3958.18  &  3657.74  \\  
 20  &  1.99  &  1846.33  &  1797.13  \\  
 30  &  1.93  &  1214.46  &  1219.85  \\  
 40  &  1.91  &  851.37  &  828.97  \\  
 50  &  1.87  &  642.29  &  618.10  \\  
 60  &  1.80  &  497.42  &  463.04  \\  
 70  &  1.77  &  368.53  &  342.28  \\  
 80  &  1.71  &  309.67  &  300.36  \\  
		\hline
	\end{tabular}
	\caption{Computational speedup  relative to the finite element method (FEM)  for the GN and FOEIM-GN methods as a function of $N$. The speedup is calculated as the ratio between the computational time of FEM and the online computational time of ROM.} 
	\label{ex3tab2}
\end{table}

\section{Conclusion}

In this work, we proposed an efficient and accurate model reduction methodology for the real-time solution of time-dependent nonlinear partial differential equations (PDEs) with parametric dependencies. Specifically, we focused on two representative examples: the Allen-Cahn equation, which governs phase separation processes, and the Buckley-Leverett equation, which models two-phase fluid flow in porous media.  We demonstrated the effectiveness of our methodology through  numerical experiments on both the Allen-Cahn and Buckley-Leverett equations. The results show that our approach delivers solutions with high accuracy while achieving computational speeds several orders of magnitude faster than the FEM. Furthermore, we explored the stability, efficiency, and accuracy of the method across different parametric variations, highlighting its robustness in a variety of application scenarios. 

While the current work focuses on two specific PDEs, the generality of our framework allows it to be extended to other complex nonlinear PDEs and different hyperreduction methods. Future work will focus on expanding the range of applications, as well as exploring adaptive strategies for further improving computational efficiency and accuracy through nonlinear manifolds. Additionally, the integration of machine learning techniques with our model reduction framework presents a promising avenue for accelerating the discovery of reduced models and optimizing hyper-reduction techniques for even more complex systems.


\section*{Acknowledgements} \label{}
I would like to thank Professors Jaime Peraire, Anthony T. Patera, and Robert M. Freund at MIT, and Professor Yvon Maday at University of Paris VI for fruitful discussions. I gratefully acknowledge a Seed Grant from the MIT Portugal Program, the United States Department of Energy under contract DE-NA0003965 and the Air Force Office of Scientific Research under Grant No. FA9550-22-1-0356 for supporting this work.  


\bibliographystyle{elsarticle-num} 
\bibliography{library.bib}




\end{document}